\documentclass[10pt]{article}

\usepackage{amsmath,amssymb,amsthm,bm}
\usepackage{geometry}
\usepackage{booktabs}
\usepackage{enumitem}
\usepackage{authblk}
\usepackage[colorlinks=true,linkcolor=blue,citecolor=blue,urlcolor=blue]{hyperref}

\newtheorem{theorem}{Theorem}[section]

\newtheorem{problem}[theorem]{Problem}
\newtheorem{remark}[theorem]{Remark}

\newcommand{\Had}{\circ}

\title{A counterexample to Huang's weak majorization conjecture}
\author[1]{Shuo Shi}
\author[1]{Juan Zhang}
\author[2]{Yun Zhang\thanks{Corresponding author. E-mail: zhangyunmaths@163.com}}
\affil[1]{School of Mathematics and Computational Science, Hunan Key Laboratory for Computation and Simulation in Science and Engineering, Key Laboratory of Intelligent Computing and Information Processing of Ministry of Education,
Xiangtan University, Xiangtan, Hunan,  411105, China}
\affil[2]{School of Mathematics and Statistics, Huaibei Normal University, Huaibei, 235000, China}
\date{}

\begin{document}

\maketitle

\begin{abstract}
We give a counterexample to a weak majorization problem, proposed by Z. Huang (Linear Algebra Appl., 434 (2) (2011) 457--462), for singular values of Hadamard products of nonnegative matrices.
The problem asks whether, for all nonnegative matrices $A$ and $B$,
\begin{equation*}
\bigl\{s_j^2(A\Had B)\bigr\}
\prec_w
\bigl\{s_j(A\Had A)s_j(B\Had B)\bigr\}
\end{equation*}
holds, where $\Had$ denotes the Hadamard product, $\prec_w$ means weak majorization, and $s_j(\cdot)$ is the $j$th largest singular value of a matrix. 
To answer this problem, we construct two $3\times3$ entrywise positive symmetric matrices. 
We give rigorous upper and lower bounds for the relevant singular-value sums by Sturm's theorem, and show that the inequality fails for the partial sum with $k=2$.


\end{abstract}

\section{Introduction}

Let $M_n$ denote the vector space of all complex $n\times n$ matrices. 
For $A\in M_n$, the singular values of $A$ are the nonnegative square roots of the eigenvalues of $A^*A$, and they are denoted by
$$
s_1(A)\geq s_2(A)\geq \cdots \geq s_n(A).
$$
For a real vector $x=(x_1,\ldots,x_n)$, let $x_{[1]}\geq x_{[2]}\geq \cdots \geq x_{[n]}$ be its components arranged in nonincreasing order. 
If $x,y\in\mathbb{R}^n$ satisfy
$$
\sum_{i=1}^{k}x_{[i]}\leq \sum_{i=1}^{k}y_{[i]},
~
k=1,2,\ldots,n,
$$
then $x$ is said to be weakly majorized by $y$, and we write $x\prec_w y$.

For two matrices $A=(a_{ij})$ and $B=(b_{ij})$ of the same size, their Hadamard product is defined by
$$
A\circ B=(a_{ij}b_{ij}).
$$
Hadamard products are closely related to several matrix inequalities involving singular values, spectral norms, spectral radii and unitarily invariant norms; 
see, for example, \cite{AndoHornJohnson1987,HornMathias1990,HornJohnson1991,Zhan1997}. 
Among such results, inequalities of Cauchy--Schwarz type for Hadamard products play a central role. 
In this direction, R. Huang \cite{Huang2008} considered the following weak majorization problem:
\begin{equation}\label{eq:problem}
\bigl\{s_j^2(A\circ B)\bigr\}
\prec_w
\bigl\{s_j(A\circ \overline{A})s_j(B\circ \overline{B})\bigr\},
~A,B\in M_n.
\end{equation}
It is known that \eqref{eq:problem} holds for the cases $k=1$ and $k=n$. 
Thus the remaining difficulty lies in the intermediate partial sums. 
This makes \eqref{eq:problem} a genuine weak majorization question rather than only a norm or trace inequality.

Z. Huang \cite{Huang2011} later showed that \eqref{eq:problem} is not valid in general for complex matrices. 
However, the case of entrywise nonnegative matrices was left open. 
More precisely, the following question was proposed in \cite[Problem~1]{Huang2011}.

\begin{problem}\label{con:Huang2011}
Does \eqref{eq:problem} hold for all nonnegative matrices $A$ and $B$?
\end{problem}

This question is natural because entrywise nonnegativity often restores inequalities that fail in the general complex setting. 
For instance, Audenaert \cite{Audenaert2010} and Horn and Zhang \cite{HornZhang2010} proved Zhan's conjecture \cite{Zhan2009} on the spectral radius of Hadamard products of nonnegative matrices. 
Therefore, it is reasonable to ask whether the weak majorization inequality \eqref{eq:problem}, although false for arbitrary complex matrices, might still be valid under entrywise nonnegativity.

The purpose of this paper is to give a negative answer to Problem \ref{con:Huang2011}. 
We construct two $3\times3$ entrywise positive symmetric matrices such that \eqref{eq:problem} fails at the intermediate partial sum $k=2$.
The verification is not merely numerical.
We compute the relevant characteristic polynomials and use Sturm's theorem to locate their roots in explicit rational intervals. 
These computations yield rigorous upper and lower bounds for the singular-value sums involved in the weak majorization inequality. 
Consequently, the proposed counterexample provides an exact disproof of Problem \ref{con:Huang2011} in the nonnegative matrix setting.

\section{The counterexample}

In the counterexample below, Sturm's theorem is used to give rigorous
bounds for the relevant singular values.

\begin{theorem}\label{lem:sturm}{\rm\cite[p. 94]{Gathen&Gerhard2013}}
Let $p(t)$ be a real polynomial. Its Sturm sequence is defined by
$$
p_0(t)=p(t),~p_1(t)=p'(t),
$$
and, for $i\geq 1$, by the Euclidean recursion
$$
p_{i-1}(t)=q_i(t)p_i(t)-p_{i+1}(t),
$$
until a constant polynomial is reached. Let $V(a)$ be the number of sign
changes in
$$
p_0(a),~p_1(a),~\ldots,~p_m(a),
$$
zeros being omitted. If neither $a$ nor $b$ is a root of $p$, then the
number of distinct real roots of $p$ in $(a,b)$ is $V(a)-V(b)$.
\end{theorem}

\begin{remark}
In the applications below, members of a Sturm sequence may be multiplied
by positive constants. This does not change the number of sign variations.
\end{remark}

We consider the following two entrywise positive symmetric matrices:
$$
A=
\begin{pmatrix}
1&4&1\\
4&3&5\\
1&5&4
\end{pmatrix},
~
B=
\begin{pmatrix}
3&4&2\\
4&4&5\\
2&5&5
\end{pmatrix}.
$$
Set
$$
C=A\Had B,~X=A\Had A,~Y=B\Had B.
$$
Then
$$
C=
\begin{pmatrix}
3&16&2\\
16&12&25\\
2&25&20
\end{pmatrix},
~
X=
\begin{pmatrix}
1&16&1\\
16&9&25\\
1&25&16
\end{pmatrix},
~
Y=
\begin{pmatrix}
9&16&4\\
16&16&25\\
4&25&25
\end{pmatrix}.
$$
We prove that
\begin{equation}\label{eq:new-violation}
s_1^2(C)+s_2^2(C)
>
s_1(X)s_1(Y)+s_2(X)s_2(Y).
\end{equation}

First,
$$
CC^T=
\begin{pmatrix}
269&290&446\\
290&1025&832\\
446&832&1029
\end{pmatrix}.
$$
Let
$$
\lambda_1\geq \lambda_2\geq \lambda_3\geq 0
$$
be the eigenvalues of $CC^T$. Since
$$
\operatorname{tr}(CC^T)=2323,
$$
we have
\begin{equation}\label{eq:new-left-trace}
s_1^2(C)+s_2^2(C)=\lambda_1+\lambda_2=2323-\lambda_3.
\end{equation}
The characteristic polynomial of $CC^T$ is
$$
p_C(t)=\det(tI-CC^T)
=t^3-2323t^2+632011t-22306729.
$$
A direct calculation gives
$$
p_C(0)=-22306729<0,
~
p_C(42)=214049>0.
$$
Hence $p_C$ has a root in $(0,42)$. Since $CC^T$ is symmetric positive
semidefinite, all roots of $p_C$ are nonnegative real eigenvalues of
$CC^T$. Therefore $\lambda_3<42$. It follows from \eqref{eq:new-left-trace}
that
\begin{equation}\label{eq:new-left-bound}
s_1^2(C)+s_2^2(C)>2323-42=2281.
\end{equation}

It remains to give a strict upper bound for the right-hand side of \eqref{eq:new-violation}. We first estimate the singular values of $X$.
We have
$$
XX^T=
\begin{pmatrix}
258&185&417\\
185&962&641\\
417&641&882
\end{pmatrix},
$$
and
$$
p_X(t)=\det(tI-XX^T)
=t^3-2102t^2+705241t-14333796.
$$
A Sturm sequence for $p_X$, normalized by positive factors, is
$$
p_X(t),
~
p_X'(t),
~
2302681t-676706209,
~
1.
$$
The signs at the required points are as follows:
$$
\begin{array}{c c c}
\toprule
t & \text{signs} & V_X(t)\\
\midrule
\frac{19533}{50} & (-,-,+,+) & 1\\
\frac{168963}{100} & (+,+,+,+) & 0\\
+\infty & (+,+,+,+) & 0\\
\bottomrule
\end{array}
$$
Here $V(+\infty)$ denotes the sign variation obtained from the leading coefficients of the corresponding Sturm polynomials.
By Theorem \ref{lem:sturm}, $p_X$ has exactly one root in
$\left(\frac{19533}{50},+\infty\right)$ and no root in
$\left(\frac{168963}{100},+\infty\right)$. Since neither
$\frac{19533}{50}$ nor $\frac{168963}{100}$ is a root of $p_X$, the
eigenvalues of $XX^T$ satisfy
$$
s_1^2(X)<\frac{168963}{100},
~
s_2^2(X)<\frac{19533}{50}.
$$
Thus
\begin{equation}\label{eq:new-X-bound}
s_1(X)<\sqrt{\frac{168963}{100}},
~
s_2(X)<\sqrt{\frac{19533}{50}}.
\end{equation}

Similarly,
$$
YY^T=
\begin{pmatrix}
353&500&536\\
500&1137&1089\\
536&1089&1266
\end{pmatrix},
$$
and
$$
p_Y(t)=\det(tI-YY^T)
=t^3-2756t^2+564484t-30041361.
$$
A Sturm sequence for $p_Y$, again normalized by positive factors, is
$$
p_Y(t),
~
p_Y'(t),
~
11804168t-1285345655,
~
1.
$$
The corresponding sign table is
$$
\begin{array}{c c c}
\toprule
t & \text{signs} & V_Y(t)\\
\midrule
\frac{2257}{20} & (-,-,+,+) & 1\\
\frac{63457}{25} & (+,+,+,+) & 0\\
+\infty & (+,+,+,+) & 0\\
\bottomrule
\end{array}
$$
Hence $p_Y$ has exactly one root in
$\left(\frac{2257}{20},+\infty\right)$ and no root in
$\left(\frac{63457}{25},+\infty\right)$. Since neither
$\frac{2257}{20}$ nor $\frac{63457}{25}$ is a root of $p_Y$, we obtain
$$
s_1^2(Y)<\frac{63457}{25},
~
s_2^2(Y)<\frac{2257}{20}.
$$
Therefore
\begin{equation}\label{eq:new-Y-bound}
s_1(Y)<\sqrt{\frac{63457}{25}},
~
s_2(Y)<\sqrt{\frac{2257}{20}}.
\end{equation}

Combining \eqref{eq:new-X-bound} and \eqref{eq:new-Y-bound} gives
$$
s_1(X)s_1(Y)+s_2(X)s_2(Y)
<
\sqrt{\frac{168963}{100}\times\frac{63457}{25}}
+
\sqrt{\frac{19533}{50}\times\frac{2257}{20}}.
$$
Moreover,
$$
\frac{168963}{100}\times\frac{63457}{25}
=
\frac{10721885091}{2500}
<
\frac{10721981209}{2500}
=
\left(\frac{103547}{50}\right)^2,
$$
and
$$
\frac{19533}{50}\times\frac{2257}{20}
=
\frac{440859810}{10000}
<
\frac{440874009}{10000}
=
\left(\frac{20997}{100}\right)^2.
$$
Thus
\begin{equation}\label{eq:new-right-bound}
s_1(X)s_1(Y)+s_2(X)s_2(Y)
<
\frac{103547}{50}+\frac{20997}{100}
=
\frac{228091}{100}.
\end{equation}

From \eqref{eq:new-left-bound} and \eqref{eq:new-right-bound}, we obtain
$$
s_1^2(C)+s_2^2(C)
>
2281
=
\frac{228100}{100}
>
\frac{228091}{100}
>
s_1(X)s_1(Y)+s_2(X)s_2(Y).
$$
Consequently,
$$
\sum_{j=1}^{2}s_j^2(A\Had B)
>
\sum_{j=1}^{2}s_j(A\Had A)s_j(B\Had B).
$$
Therefore, Problem \ref{con:Huang2011} has a negative answer even for entrywise positive symmetric matrices.

\section{Conclusion}

In this paper, we have given a counterexample to the weak majorization problem for singular values of Hadamard products of nonnegative matrices proposed by Z. Huang \cite[Problem~1]{Huang2011}. 
This counterexample also indicates that the validity of Cauchy--Schwarz type inequalities for Hadamard products is sensitive to the form of the singular-value comparison. 
Although the first and the full partial sums are known to satisfy the proposed inequality, the present example shows that the intermediate partial sums may behave differently. 
Therefore, any positive result in this direction would require additional structural assumptions on the matrices, or a modification of the comparison term on the right-hand side.

~\\

\noindent{\bf Acknowledgement.} The work was supported by the National Science Foundation of Anhui Higher Education Institutions of China (KJ2021ZD0058).

\end{document}